\documentclass[12pt,a4paper]{amsart}
\usepackage{amssymb,amsmath}
\usepackage{mathrsfs}

\headsep=1cm \numberwithin{equation}{section}
\newtheorem{theorem}{Theorem}[section]
\newtheorem{lemma}[theorem]{Lemma}
\newtheorem{proposition}[theorem]{Proposition}
\newtheorem{corollary}[theorem]{Corollary}

\theoremstyle{definition}

\newtheorem{example}[theorem]{Example}

\newcommand\Ass{\operatorname{Ass}}
\newcommand\mAss{\operatorname{mAss}}

\newcommand\Spec{\operatorname{Spec}}
\newcommand\Rad{\operatorname{Rad}}

\newcommand\height{\operatorname{height}}

\begin{document}

\title[A Characterization of locally quasi-unmixed rings]{A Characterization of locally quasi-unmixed rings}%
\author{Simin Mollamahmoudi, Adeleh Azari  and Reza Naghipour$^*$}%
\address{Department of Mathematics, University of Tabriz,
Tabriz, Iran, and, School of Mathematics, Institute for Research in Fundamental
Sciences (IPM), P.O. Box: 19395-5746, Tehran, Iran.}%
\email{mahmoudi.simin@yahoo.com (Simin Mollamahmoudi)}
\email{adeleh\_azari@yahoo.com (Adeleh Azari)}
\email{naghipour@ipm.ir (Reza Naghipour)} \email {naghipour@tabrizu.ac.ir (Reza Naghipour)}

\thanks{ 2010 {\it Mathematics Subject Classification}: 13D45, 14B15, 13E05.\\
This research was in part supported by a grant from IPM. \\
$^*$Corresponding author: e-mail: {\it naghipour@ipm.ir} (Reza
Naghipour)}%
\keywords{Associated primes, ideal topologies, integral closure, locally quasi-unmixed ring, Rees ring.}

\begin{abstract}
Let $\bar{I}$ denote the integral closure of an ideal in a  Noetherian ring $R$. The main result of this paper asserts that $R$  is locally quasi-unmixed if and only if, the topologies defined by $\overline{I^n}$  and $I^{\langle n\rangle}$, $\ n\geq 1$,  are equivalent. In addition, some results about the behavior of linearly equivalent  topologies of ideals under various ring homomorphisms are included.

\end{abstract}
\maketitle
\section {Introduction}
Let $R$ denote a commutative Noetherian ring, and $I$ an ideal of $R$. The interesting concept of quintasymptotic prime ideals of $I$ was introduced by McAdam \cite{Mc2}. A prime ideal $ {\mathfrak{p}}$ of $R$ is called a {\it quintasymptotic prime ideal of} $I$ if there exists $z\in\mAss _{R_\mathfrak{p}^*} R_\mathfrak{p}^*$ such that $\Rad(I{R_{\mathfrak{p}}^*}+z)={\mathfrak{p}}R_{\mathfrak{p}}^*$. The set of  quintasymptotic prime ideals  of $I$ is denoted by $\bar{Q^*}(I)$, and it is a finite set. Also, in \cite{Ra2}, Ratliff, Jr., introduced set of associated primes $\bar{A^*} (I):= \Ass_R R/\overline{I^{n}}$ for large $n$, called the {\it presistent prime ideals} of $I$, and he showed that this  finite set has some nice properties in the theory of asymptotic prime divisors; here for any ideal $J$ of $R$,  $\bar{J}$ denotes   the {\it integral closure}  of $J$ in  $R$, i.e.,  $\bar{J}$  is the ideal of $R$ consisting of all elements $x\in R$ which satisfy an equation $x^n+r_1x^{n-1}+\dots+r_n=0$, where $r_i\in J^i,\ i=1,\ldots,n$.

In his famous paper \cite{Re}, D. Rees showed  that a local ring $(R,\mathfrak{m})$ is analytically unramified  if and only if the topology defined by $\overline{I^{n}},\ n\geqslant1$, is equivalent to the  $I$-adic topology for an $\mathfrak{m}$-primary ideal  $I$ of   $R$. In  \cite{Ra1}, L.J. Ratliff, Jr., proved corresponding results in order to characterize reduced unmixed local rings. (Recall that a local ring $(R,\mathfrak{m})$ is  called {\it analytically unramified} (resp. {\it unmixed}), if the $\mathfrak{m}$-adic completion, $R^*$,  of $R$ is reduced (resp. all the prime ideals of $\Ass_{R^*}R^*$ have the same dimension)). The main theorem of this paper gives a characterization of locally quasi-unmixed  Noetherian rings, which is closely related to Ress' result \cite{Re}. Since such rings occur in many investigations in commutative algebra and algebraic geometry, it is desirable to know as many  properties of such ring as possible. This characterization gives one such property, and that such rings have this property is a new result, and until now was not know to hold even in a regular local ring. More precisely we shall show that:
\begin{theorem}\label{thm1}
Let $R$ denote a commutative Noetherian ring. Then the following conditions are equivalent:

\begin{itemize}
\item[(i)]  $R$  is locally quasi-unmixed.\\\vspace{-.3cm}
\item[(ii)] For every ideal $I$ of the principal class in $R$,  the   topologies defined   by  $\overline{I^{n}}$ and $I^{\langle n\rangle}$, $n\geq 1$, are linearly equivalent.
\item[(iii)]
 For every ideal $I$ of the principal class in $R$,  the   topologies defined   by  $\overline{I^{n}}$ and $I^{\langle n\rangle}$, $n\geq 1$, are equivalent.
\end{itemize}
\end{theorem}
 This is closely related to Ress' result \cite{Re} for  characterization  quasi-unmixed local  rings.  Here $I^{\langle n\rangle}$ denotes the union $(\overline{I^{n}}:_Rs)$, where $s$ varies in $R\backslash \bigcup \{\frak p\in \mAss_RR/I\}$. See Theorem \ref{thm1} for the proof of Theorem 1.1.

 One of our tools for proving Theorem 1.1 is the following, which is a  characterization of the equivalence between the topologies defined  by the filtration $\overline{I^{n}}$ and $I^{\langle n\rangle}$, $n\geqslant1$.
 
\begin{proposition}\label{prop1}
  Let  $I$ denote  an ideal in a commutative Noetherian ring $R$.  Then the topologies  defined by $\overline{I^{n}}$ and $I^{\langle n\rangle}$, $n\geq 1$, are  equivalent {\rm(}resp. linearly equivalent{\rm)} if and only if $\bar{Q^*}(I)$ {\rm(}resp. $\bar{A^*}(I)${\rm)} is equal to $\mAss _R R/I$.
\end{proposition}

Pursuing this point of view further we prove some results about the behavior of the linearly equivalent topologies of ideals under various ring homomorphisms. In connection to this we derive the following consequence of Proposition \ref{prop1}.

\begin{corollary}\label{coro1}
  Let   $R$ be a  Noetherian ring and $T$ be a  finitely generated integral ring extension of $R$   such that every minimal prime of $T$ lies over a  minimal prime of $R$.  If the topologies $\overline{I^n}$ and $I^{\langle n\rangle},\ n\geqslant1,$ are  linearly  equivalent,   then the topologies defined by $\overline{(IT)^n}$ and $(IT)^{\langle n\rangle},\  n\geqslant1,$ are also linearly equivalent; and the converse holds whenever $T$ is faithfully flat.
\end{corollary}

  Throughout  this paper, for any commutative Noetherian ring $R$ with nonzero identity, and  for any ideal $I$ of $R$, we denote by $\mAss_R R/I$ the set of minimal prime  ideals over $I$. If $(R,\mathfrak{m})$ is local, then  $R^* $ denotes the completion of $R$ with respect to the $\mathfrak{m}$-adic topology. Then $R$ is said to be {\it quasi-unmixed ring} if  for  every $\mathfrak{p}\in \mAss_{R^*}R^*$, the condition $\dim R^*/\mathfrak{p}=\dim R$ is satisfied. More generally, if $R$ is not necessarily local, $R$ is a {\it locally quasi-unmixed ring}  if for any  $\mathfrak{p}\in\Spec(R)$, $R_{\mathfrak{p}}$ is a local quasi-unmixed ring. For any ideal $I$ of $R$, we denote by $\mathscr{R}$ the graded Rees ring $R[u,It]:= \bigoplus\limits_{n\in \mathbb{Z}}I^nt^n$ of $R$ with respect to $I$, where $t$ is an indeterminate and $u=t^{-1}$.
  Also, the {\it radical of} $I$, denoted by $\Rad(I)$, is defined to be the set $\{x\in R \mid x^n\in I\ \textrm{for some}\ n\in \mathbb{N}\}$. Finally, if $(R,\mathfrak{m})$ is local, then the {\it analytic spread of} $I$ is defined to be $\ell(I):=\dim \mathscr{R}/(\mathfrak{m},u)\mathscr{R}$ (see \cite{NR}).  For any unexplained notation and terminology we refer the reader to \cite{Ma} or \cite{Nag}.\\

 \section{Locally quasi-unmixed rings and comparison of topologies  }\label{sec1}
The purpose of this section is to establish a characterization of locally quasi-unmixed  Noetherian rings, which is closely related to Ress' result \cite{Re}. The main goal is Theorem \ref{thm1}. The following lemmas are needed in the proof of that theorem.
\begin{lemma}\label{lem1}
Let $I$ be an ideal of  a Noetherian ring $R$. Then the following conditions are equivalent:
\begin{itemize}
  \item[(i)]  $\bar{Q^*}(I)=\mAss_RR/I$.\\\vspace{-.3cm}
  \item[(ii)] The   topologies defined   by  $\overline{I^{n}}$ and $I^{\langle n\rangle}$, $n\geq 1$, are  equivalent.
\end{itemize}
\end{lemma}
\proof  The assertion follows easily from \cite[Theorem 1.5]{Mc2} and the fact that $\mAss_RR/I\subseteq \bar{Q^*}(I)$.   \qed

\begin{lemma}\label{lem2}
  Let   $R$ be a  Noetherian ring such that   the topologies  defined   by $\overline{I^{n}}$ and $I^{\langle n\rangle},\ n\geqslant1,$ are    equivalent   for all ideals $I$ of the principal class in $R$. Then for every prime ideal $\mathfrak{p}$ of $R$ and every ideal $J$ of the principal class in $R_{\mathfrak{p}}$, the topologies  defined   by $\overline{J^{n}}$ and $J^{\langle n\rangle},\ n\geqslant1,$ are equivalent.
  \end{lemma}
  \proof  Let $\mathfrak{p}\in\Spec(R)$ and let $J$ be an ideal of the principal class in $R_{\mathfrak{p}}$. Then in view of \cite[Lemma 5.1]{V1}, there exists an ideal $I$ of $R$ of the principal class such that $J=IR_{\mathfrak{p}}$. Now, in view of Lemma \ref{lem1}, it is enough for us to show that  $\bar{Q^*}(J)=\mAss_{R_{\mathfrak{p}}}R_{\mathfrak{p}}/J$. To do this, let $\mathfrak{q}R_{\mathfrak{p}}\in  \bar{Q^*}(J)$. That is   $\mathfrak{q}R_{\mathfrak{p}}\in  \bar{Q^*}(IR_{\mathfrak{p}})$. Then, by \cite[Proposition 1.1]{Mc2}, $\mathfrak{q}\in  \bar{Q^*}(I)$, and so by Lemma \ref{lem1}, $\mathfrak{q}\in \mAss_RR/I$. Therefore, $\mathfrak{q}R_{\mathfrak{p}}\in \mAss_{R_{\mathfrak{p}}}R_{\mathfrak{p}}/IR_{\mathfrak{p}}$, as required.  \qed\\
  
The next Lemma was proved by McAdam and Ratliff in \cite{MR1}.

\begin{lemma}\label{lem3}
Let $I$ be an ideal of  a locally quasi-unmixed Noetherian ring $R$ such that $\ell(IR_{\mathfrak{p}})=\height(IR_{\mathfrak{p}})$ for all $\mathfrak{p}\in \bar{A^*}(I)$. Then  $\bar{A^*}(I)=\mAss_RR/I$.
\end{lemma}
\proof See \cite[Lemma 5.4]{MR1}.  \qed\\

\begin{lemma}\label{lem4}
Let $I$ denote  an ideal in a Noetherian ring $R$.  Then the following conditions are equivalent:
\begin{itemize}
  \item[(i)]  $\bar{A^*}(I)=\mAss_RR/I$.\\\vspace{-.3cm}
  \item[(ii)] The   topologies defined   by  $\overline{I^{n}}$ and $I^{\langle n\rangle}$, $n\geq 1$, are linearly equivalent.
\end{itemize}
\end{lemma}
\proof The result follows from \cite[Corollary 1.6]{Mc2} and the fact that  $\mAss_RR/I\subseteq \bar{A^*}(I)$. \qed\\

We are now ready to state and prove the main theorem of this section which is a characterization of   locally quasi-unmixed Noetherian rings in terms of the equivalence (resp. linearly equivalence)  between the topologies induced by $\overline{I^{n}}$ and $I^{\langle n\rangle}$, $n\geq 1$, for the principal class ideals $I$ of $R$. Recall that an ideal $I$ of $R$ is called {\it of the principal class} if $I$ is generated by $\height I$ elements.

\begin{theorem}\label{thm1}
Let $R$ denote a commutative Noetherian ring. Then the following conditions are equivalent:
\begin{itemize}
\item[(i)]  $R$  is locally quasi-unmixed.\\\vspace{-.3cm}
\item[(ii)] For every ideal $I$ of the principal class in $R$,  the   topologies defined   by  $\overline{I^{n}}$ and $I^{\langle n\rangle}$, $n\geq 1$, are linearly equivalent.
\item[(iii)]
 For every ideal $I$ of the principal class in $R$,  the   topologies defined   by  $\overline{I^{n}}$ and $I^{\langle n\rangle}$, $n\geq 1$, are equivalent.
\end{itemize}
\end{theorem}
\proof
 First we show $\rm(i)\Longrightarrow (ii)$. If $R$ is locally quasi-unmixed, then in view of Lemmas \ref{lem3} and  \ref{lem4}, it is  enough for us to show that, for all $\mathfrak{p}\in\bar{A^*}(I)$, $\ell(IR_{\mathfrak{p}})=\height(IR_{\mathfrak{p}})$ for every ideal $I$ of the principal class. To this end, in view of \cite[Proposition 4.1]{Mc1}, $\height({\mathfrak{p}})=\ell(IR_{\mathfrak{p}})$. Now, since at least $\ell(\mathfrak{a})$ elements are needed to generate $\mathfrak{a}$, for any ideal $\mathfrak{a}$ in a commutative Noetherian ring $A$, and as $IR_{\mathfrak{p}}$ is an ideal of the principal class in  $R_{\mathfrak{p}}$, it follows that $\ell(IR_{\mathfrak{p}})\leq\height(IR_{\mathfrak{p}})$. Furthermore, since $I\subseteq \mathfrak{p}$, it yields that   $$\height(IR_{\mathfrak{p}})\leq\height(\mathfrak{p}R_{\mathfrak{p}})=\height({\mathfrak{p}}),$$ and so $$\ell(IR_{\mathfrak{p}})\leq\height(IR_{\mathfrak{p}})\leq\height({\mathfrak{p}})=\height(IR_{\mathfrak{p}}).$$
 Hence $\ell(IR_{\mathfrak{p}})=\height(IR_{\mathfrak{p}})$, as required.

 Now, because of the implication $\rm(ii)\Longrightarrow(iii)$ is trivially true, so in order to complete the proof we have to show that  $\rm(iii)\Longrightarrow(i)$.
 Let $\mathfrak{p}\in\Spec(R)$. We need to show that $R_{\mathfrak{p}}$ is a quasi-unmixed ring. To do this, in view of Lemma \ref{lem2} and \cite[Remark 2.9]{Ra}, without loss of generality we may assume that $(R,\mathfrak{m})$ is a local ring. Now, for proving the quasi-unmixedness of $R$, there are two cases to consider.

 \emph{\bf{Case 1.}} Suppose that  $\mathfrak{m}R^*\in \mAss_{R^*}R^*$. Then $\height(\mathfrak{m}R^*)=0$, and so $\dim R^*=0$. Hence, $R$ is a quasi-unmixed ring, as required.

  \emph{\bf{Case 2.}} Now, suppose that  $\mathfrak{m}R^*\notin \mAss_{R^*}R^*$, and let $\mathfrak{q}\in \mAss_{R^*}R^*$. We need to show that $\dim R^*/\mathfrak{q}=\dim R$. To this end, as  $\mathfrak{m}R^*\notin \mAss_{R^*}R^*$, we have $\dim R^*/\mathfrak{q}:=n$, where $n>0$. Therefore in view of \cite[Proposition 3.5]{Na1}, there exists an ideal $\mathfrak{a}$ of $R$  of the principal class of height $n$ and $\Rad(\mathfrak{a}R^*+\frak q)=\mathfrak{m}R^*$. Whence, $\mathfrak{m}\in \bar{Q^*}(\mathfrak{a})$. Moreover, as $\mathfrak{a}$ is  the principal class, it follows from assumption $\rm(iii)$ and Lemma \ref{lem1} that $\mathfrak{m}\in \mAss _R R/\mathfrak{a}$. Consequently, $\height({\mathfrak{m}})=n$, and so $\dim R^*/\mathfrak{q}=\dim R$, as required. \qed\\

  The following corollary gives us a  characterization of   locally quasi-unmixed Noetherian rings in terms of  quintasymptotic and  presistent prime ideals of $I$.
\begin{corollary}\label{coro2}
Let $R$ be a commutative  Noetherian ring. Then the following conditions are equivalent:
\begin{itemize}
  \item[(i)]  $R$ is locally quasi-unmixed.\\\vspace{-.3cm}
  \item[(ii)]  $\bar{A^*}(I)=\mAss_RR/I$, for every ideal $I$ of the principal class of  $R$.
  \item[(iii)] $\bar{Q^*}(I)=\mAss_RR/I$, for every ideal $I$ of the principal class of  $R$.
\end{itemize}
\end{corollary}
\proof The assertion follows from Theorem \ref{thm1}, Lemma \ref{lem1} and \cite[Lemma 2.1]{Mc2}. \qed\\

As the final result of this section, we construct an example to show that the Theorem \ref{thm1} is not true, if $I$  is not ideal of the  principal  class. The following lemma is needed in the proof of the Example \ref{exam1}.
\begin{lemma}\label{lem5}
Let $R$ be a Noetherian ring such that $\dim R>0$.  Let $I\subseteq \mathfrak{p}$ be ideals of $R$ with $\mathfrak{p}\in\Spec(R)$. Then  the following conditions are equivalent:
\begin{itemize}
  \item[(i)]  $\mathfrak{p}\in\bar{A^*}(I)$.\\\vspace{-.3cm}
  \item[(ii)] $\mathfrak{p}\in\bar{A^*}(xI)$ for any element $x$ not contained in any minimal prime of $R$.
\end{itemize}
\end{lemma}
\proof  See \cite[Proposition 3.26]{Mc1}. \qed\\
\begin{example}\label{exam1}
Let $k$ be a field and let $R=k[x,y]_{(x,y)}$. Set $\mathfrak{m}=(x,y)R$ and $I=x\mathfrak{m}$. Then $\mathfrak{m}\in \bar{A^*}(I)$ and $\mathfrak{m}\notin \bar{Q^*}(I)$.
\end{example}
\proof Since $x$ is not contained in any minimal prime of $R$ and $\mathfrak{m}\in \bar{A^*}(\mathfrak{m})$, it follows from Lemma \ref{lem5} that  $\mathfrak{m}\in \bar{A^*}(I)$. Now, we need to show that  $\mathfrak{m}\notin \bar{Q^*}(I)$. Suppose, the contrary, that $\mathfrak{m}\in \bar{Q^*}(I)$. Then, there exists $z\in \mAss_{R^*}R^*$ such that $\Rad(IR^*+z)=\mathfrak{m}R^*$. Since $I=x\mathfrak{m}$, it yields that $\Rad(xR^*+z)=\mathfrak{m}R^*$. Hence, $\mathfrak{m}R^*/z$ is minimal over $x(R^*/z)$, and so in view of Krull's Principal Ideal Theorem,
$\height (\mathfrak{m}R^*/z)\leq1$. On the other hand, as $R^*$ is a Cohen-Macaulay ring, it follows that
$$\height (\mathfrak{m}R^*/z)=\height (\mathfrak{m}R^*)-\height (z),$$
and so $\height (\mathfrak{m}R^*/z)=2$, which provides a contradiction.  \qed\\

\section{Lineally equivalent topologies}
 Our aim of this section is to obtain some results about the behavior of the lineally equivalent topologies of ideals under various ring homomorphisms.\\
 
\begin{proposition}\label{prop2}
Let $R$ be a Noetherian ring and let $I$ be  an ideal of $R$ such that  the topologies  defined by $\overline{I^{n}}$ and $I^{\langle n\rangle},\ n\geq 1,$ are linearly equivalent. Let $T$ be a  finitely generated integral ring extension of $R$ such that every minimal prime of $T$ lies over a  minimal prime of $R$. Then the topologies  defined by $\overline{(IT)^n}$ and $(IT)^{\langle n\rangle}$, $n\geq1$,  are linearly equivalent.
\end{proposition}
\proof In view of Lemma \ref{lem4}, it is enough to show that $\bar{A^*}(IT)=\mAss_TT/IT$. To this end, let $\mathfrak{p}\in\bar{A^*}(IT)$ and we show that $\mathfrak{p}\in\mAss_TT/IT$. Suppose the contrary that $\mathfrak{p}\notin\mAss_TT/IT$. Then, there exists $\mathfrak{q}\in\mAss_TT/IT$ such that $\mathfrak{q}\subsetneqq\mathfrak{p}$. Now as $\mathfrak{p}, \mathfrak{q}\in\bar{A^*}(IT)$, ( note that $\mAss_TT/IT\subseteq \bar{A^*}(IT)$), it follows from \cite[Theorem 3.3]{Ra2} that $\mathfrak{p}\cap R$ and $\mathfrak{q}\cap R$ are  contained in $\bar{A^*}(I)$. Hence, in view of Lemma \ref{lem4},  $\mathfrak{p}\cap R$ and $\mathfrak{q}\cap R$ are  contained in $\mAss_RR/I$, and so $\mathfrak{p}\cap R=\mathfrak{q}\cap R$. Therefore, by \cite[Theorem 9.3]{Ma}, $\mathfrak{p}=\mathfrak{q}$ which is a contradiction.  \qed\\
\begin{proposition}\label{prop3}
Let $R$ be a Noetherian ring and let $I$ be  an ideal of $R$. Let $T$ be a faithfully flat ring extension of $R$ such that the topologies  defined by $\overline{(IT)^n}$ and $(IT)^{\langle n\rangle}, n\geq 1,$ are linearly equivalent. Then the topologies  defined by $\overline{I^{n}}$ and $I^{\langle n\rangle}$, $n\geq1$,  are linearly equivalent.
\end{proposition}
\proof In view of Lemma \ref{lem4}, it is enough to show that $\bar{A^*}(I)=\mAss_RR/I$. To do this, let $\mathfrak{p}\in\bar{A^*}(I)$. Then in view of \cite[Corollary 6.9]{Ra1}, there exists $\mathfrak{p}^*\in\bar{A^*}(IT)$ such that $\mathfrak{p}^*\cap R=\mathfrak{p}$.  Now by virtue of  Lemma \ref{lem4}, $\mathfrak{p}^*\in\mAss_TT/IT$. Hence, by the Going Down property between $T$ and $R$ (cf. \cite[Theorem 9.5]{Ma}), we see that  $\mathfrak{p}\in\mAss_RR/I$, as required.  \qed\\
\begin{theorem}\label{thm2}
Let $R$ be a Noetherian ring and let $I$ be  an ideal of $R$.  Then the topologies  defined by $\overline{I^{n}}$ and $I^{\langle n\rangle}, n\geq 1$, are linearly equivalent if only if the topologies  defined by $\overline{(IR[x])^n}$ and $(IR[x])^{\langle n\rangle}, n\geq 1$, are linearly equivalent
\end{theorem}
\proof Since $R[x]$ is a faithfully flat ring extension of $R$, the sufficiency follows from Proposition  \ref{prop2}. For necessity,  in view of Lemma \ref{lem4},  it is enough to show that $$\bar{A^*}(IR[x])=\mAss_{R[x]}R[x]/IR[x].$$ To this end, let $\mathfrak{p}R[x]\in\bar{A^*}(IR[x])$, note that by \cite[Proposition 3.21]{Mc1},
$$\bar{A^*}(IR[x])=\{\mathfrak{p}R[x] \mid \mathfrak{p}\in \bar{A^*}(I)\}.$$
 Then $\mathfrak{p}\in \bar{A^*}(I)$, and so by  Lemma \ref{lem4},  $\mathfrak{p}\in\mAss_RR/I$. Now, it easy to see that $\mathfrak{p}R[x]\in\mAss_{R[x]}R[x]/IR[x]$, as required. \qed\\
\begin{proposition}\label{prop4}
Let $R$ be a Noetherian ring and let $I$ be  an ideal of $R$ such that the topologies  defined by $\overline{I^{n}}$ and $I^{\langle n\rangle}, n\geq 1,$ are linearly equivalent. Then for any $z\in\mAss_RR/I$, the topologies  defined by $\overline{(I+z/z)^n}$ and $(I+z/z)^{\langle n\rangle}, n\geq 1,$ are linearly equivalent
\end{proposition}
\proof In view of Lemma \ref{lem4}, it  suffices to show that $$\bar{A^*}(I+z/z)=\mAss_{R/z}((R/z)/(I+z/z)).$$ To do this, let $\mathfrak{p}/z\in\bar{A^*}(I+z/z)$. Then in view of \cite[Corollary 6.3]{Ra1},  $\mathfrak{p}\in\bar{A^*}(I)$.  Hence, by  Lemma \ref{lem4}, $\mathfrak{p}\in\mAss_RR/I$.
Now, it easy to see that $$\mathfrak{p}/z\in\mAss_{R/z}((R/z)/(I+z/z)),$$ as required. \qed\\

\begin{center}
{\bf Acknowledgments}
\end{center}
The authors would like to thank Professor Monireh Sedghi for reading of the first draft and valuable discussions.
Finally, the authors would like to thank from the Institute for Research in Fundamental Sciences (IPM), for the financial support.

\end{document}